\documentclass[12pt]{article}
\usepackage{amssymb,amsmath}
\usepackage{hyperref}
\usepackage{graphicx}
\usepackage{color}

\begin{document}

\title{\LARGE\bf A formula for pi involving nested radicals}

\author{
\normalsize\bf S. M. Abrarov\footnote{\scriptsize{Dept. Earth and Space Science and Engineering, York University, Toronto, Canada, M3J 1P3.}}\, and B. M. Quine$^{*}$\footnote{\scriptsize{Dept. Physics and Astronomy, York University, Toronto, Canada, M3J 1P3.}}}

\date{April 17, 2018}
\maketitle
\vspace{-0.75cm}

\begin{abstract}
\vspace{0.25cm}
\normalsize
We present a new formula for pi involving nested radicals with rapid convergence. This formula is based on the arctangent function identity with argument \footnotesize $x=\sqrt{2-{{a}_{k-1}}}/{{a}_{k}}$ \normalsize, where
\footnotesize
\[
{{a}_{k}}=\underbrace{\sqrt{2+\sqrt{2+\sqrt{2+\cdots +\sqrt{2}}}}}_{k\,\,\text{square}\,\,\text{roots}}
\]
\normalsize
is a nested radical consisting of $k$ square roots. The computational test we performed reveals that the proposed formula for pi provides a significant improvement in accuracy as the integer $k$ increases.
\\
\\
\noindent{\bf{Keywords:}} constant pi, arctangent function, nested radical
\\
\noindent{\bf{MSC classes:}} 11Y60

\end{abstract}

\section{Introduction}

In $1593$ the French mathematician Fran\c{c}ois Vi\'{e}te discovered a classical formula for the constant pi that can be expressed elegantly in nested radicals consisting of square roots of twos \cite{Servi2003, Levin2005, Kreminski2008}
\begin{equation}\label{eq_1}
\frac{2}{\pi }=\frac{\sqrt{2}}{2}\frac{\sqrt{2+\sqrt{2}}}{2}\frac{\sqrt{2+\sqrt{2+\sqrt{2}}}}{2}\cdots.
\end{equation}
He found this formula for pi geometrically by considering a regular polygon enclosed inside the circle with unit radius. It is convenient to define the nested radicals by recurrence relations
$$
{{a}_{k}}=\sqrt{2+{{a}_{k-1}}}\,\,\text{and}\,\,{{a}_{1}}=\sqrt{2}
$$
to represent the Vi\'{e}te\text{'}s formula \eqref{eq_1} for pi in a more compact form as
$$
\frac{\pi }{2}=\underset{K\to \infty }{\mathop{\lim }}\,\prod\limits_{k=1}^{K}{\frac{{{a}_{k}}}{2}}.
$$

Many interesting formulas involving nested radicals consisting of square roots of twos have been reported in the modern literature (see for example \cite{Servi2003, Levin2005, Kreminski2008}). In this paper we derive a new formula for pi that is also based on this type of the nested radicals
\small
\begin{equation}\label{eq_2}
\begin{aligned}
\pi = &\,\, {{2}^{k+1}} i\,\sum\limits_{m=1}^{\infty}\frac{1}{2m-1}\\
&\times\left( \frac{1}{{{\left( 1+2i\,\,{{a}_{k}}/\sqrt{2-{{a}_{k-1}}} \right)}^{2m-1}}}-\frac{1}{{{\left( 1-2i\,\,{{a}_{k}}/\sqrt{2-{{a}_{k-1}}} \right)}^{2m-1}}} \right).
\end{aligned}
\end{equation}
\normalsize
The computational test reveals that accuracy of the constant pi in the truncated formula \eqref{eq_2} can be considerably improved by successive increment of the integer $k$.

\section{Derivation}

As it has been shown in our recent paper \cite{Abrarov2016a}, any function $f\left( t \right)$ differentiable within the interval $t\in \left[ 0,1 \right]$ can be integrated numerically by truncating the parameters $L$ or $M$ in the following limits (see equations (9) and (6) in \cite{Abrarov2016a})
\[
\tag{3a}\label{eq_3a}
\int\limits_{0}^{1}{f\left( t \right)dt}=\underset{L\to \infty }{\mathop{\lim }}\,\sum\limits_{\ell =1}^{L}{\sum\limits_{m=0}^{M}{\frac{{{\left( -1 \right)}^{m}}+1}{{{\left( 2L \right)}^{m+1}}\left( m+1 \right)!}}}{{\left. {{f}^{\left( m \right)}}\left( t \right) \right|}_{t=\frac{\ell -1/2}{L}}}
\]
and
\[
\tag{3b}\label{eq_3b}
\int\limits_{0}^{1}{f\left( t \right)dt}=\underset{M\to \infty }{\mathop{\lim }}\,\sum\limits_{\ell =1}^{L}{\sum\limits_{m=0}^{M}{\frac{{{\left( -1 \right)}^{m}}+1}{{{\left( 2L \right)}^{m+1}}\left( m+1 \right)!}}}{{\left. {{f}^{\left( m \right)}}\left( t \right) \right|}_{t=\frac{\ell -1/2}{L}}},
\]
respectively.

Since
$$
\int\limits_{0}^{1}{\frac{x}{1+{{x}^{2}}{{t}^{2}}}dt}=\arctan \left( x \right)
$$
substituting the integrand
$$
f\left( t \right)=\frac{x}{1+{{x}^{2}}{{t}^{2}}}
$$
into the equation \eqref{eq_3a} results in (see \cite{Abrarov2016b} for more details in derivation)
\footnotesize
\[
\tag{4a}\label{eq_4a}
\begin{aligned}
&\hspace{-0.35cm}\arctan \left( x \right)=\\
& \,\, i\underset{L\to \infty }{\mathop{\lim }}\,\sum\limits_{\ell =1}^{L}{\sum\limits_{m=1}^{\left\lfloor \frac{M}{2} \right\rfloor +1}{\frac{1}{2m-1}\left( \frac{1}{{{\left( \left( 2\ell -1 \right)+2iL/x \right)}^{2m-1}}}-\frac{1}{{{\left( \left( 2\ell -1 \right)-2iL/x \right)}^{2m-1}}} \right)}}.
\end{aligned}
\]
\normalsize
According to equations \eqref{eq_3a} and \eqref{eq_3b} the parameters $L$ and $M$ under the limit notation are interchangeable (compare equations (9) and (6) from the paper \cite{Abrarov2016a}). Consequently, using absolutely same derivation procedure as described in \cite{Abrarov2016b} for the equation \eqref{eq_4a} above, we can also write
\footnotesize
\[
\tag{4b}\label{eq_4b}
\begin{aligned}
&\hspace{-0.3cm}\arctan\left( x \right)=\\
& \,\, i\underset{M\to \infty }{\mathop{\lim }}\,\sum\limits_{\ell =1}^{L}{\sum\limits_{m=1}^{\left\lfloor \frac{M}{2} \right\rfloor +1}{\frac{1}{2m-1}\left( \frac{1}{{{\left( \left( 2\ell -1 \right)+2iL/x \right)}^{2m-1}}}-\frac{1}{{{\left( \left( 2\ell -1 \right)-2iL/x \right)}^{2m-1}}} \right)}}.
\end{aligned}
\]
\normalsize

Comparing equations \eqref{eq_4a} and \eqref{eq_4b} we can see that at least one of the parameters $L$ or $M$ must be large enough in truncation for high-accuracy approximation. However, as the parameter $M$ is more important for rapid convergence, the limit \eqref{eq_4b} is preferable for numerical analysis. Since in the limit \eqref{eq_4b} the integer $L$ may not be necessarily large, in order to simplify it we can choose any small value, say $L=1$. This leads to the equation
\setcounter{equation}{4}
\small
\[
\arctan \left( x \right)=i\underset{M\to \infty }{\mathop{\lim }}\,\sum\limits_{m=1}^{\left\lfloor \frac{M}{2} \right\rfloor +1}{\frac{1}{2m-1}\left( \frac{1}{{{\left( 1+2i/x \right)}^{2m-1}}}-\frac{1}{{{\left( 1-2i/x \right)}^{2m-1}}} \right)}.
\]
\normalsize
Since the integer $M$ tends to infinity, the upper bound $\left\lfloor M/2 \right\rfloor +1$ of summation also tends to infinity. Consequently, this equation can be rewritten in a simplified form
\small
\begin{equation}\label{eq_5}
\arctan \left( x \right)=i\,\sum\limits_{m=1}^{\infty}{\frac{1}{2m-1}\left( \frac{1}{{{\left( 1+2i/x \right)}^{2m-1}}}-\frac{1}{{{\left( 1-2i/x \right)}^{2m-1}}} \right)}.
\end{equation}
\normalsize

Using the identity for the cosine double-angle
$$
\cos \left( \frac{\pi}{{{2}^{k}}} \right)=2{{\cos }^{2}}\left( \frac{\pi}{{{2}^{k+1}}} \right)-1
$$
and taking into consideration that
$$
\cos \left( \frac{\pi}{{{2}^{2}}} \right)=\frac{\sqrt{2}}{2},
$$
we can readily find by induction
$$
\cos \left( \frac{\pi}{{{2}^{3}}} \right)=\frac{\sqrt{2+\sqrt{2}}}{2},
$$
$$
\cos \left( \frac{\pi}{{{2}^{4}}} \right)=\frac{\sqrt{2+\sqrt{2+\sqrt{2}}}}{2},
$$
$$
\cos \left( \frac{\pi}{{{2}^{5}}} \right)=\frac{\sqrt{2+\sqrt{2+\sqrt{2+\sqrt{2}}}}}{2},
$$
$$
\vdots
$$
\begin{equation}\label{eq_6}
\cos \left( \frac{\pi}{{{2}^{k+1}}} \right)=\frac{1}{2}\underbrace{\sqrt{2+\sqrt{2+\sqrt{2+\cdots +\sqrt{2}}}}}_{k\,\,\text{square}\,\,\text{roots}}.
\end{equation}

Consider the following relation
\begin{equation}\label{eq_7}
\frac{\pi }{{{2}^{k+1}}}=\arctan \left( \tan \left( \frac{\pi}{{{2}^{k+1}}} \right) \right)=\arctan \left( \frac{\sin \left( \frac{\pi}{{{2}^{k+1}}} \right)}{\cos \left( \frac{\pi}{{{2}^{k+1}}} \right)} \right),
\end{equation}
where
\begin{equation}\label{eq_8}
\sin \left( \frac{\pi}{{{2}^{k+1}}} \right)=\sqrt{1-{{\cos }^{2}}\left( \frac{\pi}{{{2}^{k+1}}} \right)}.
\end{equation}
Substituting the equation \eqref{eq_6} into the identity \eqref{eq_8} leads to
\begin{equation}\label{eq_9}
\sin \left( \frac{\pi}{{{2}^{k+1}}} \right)=\frac{1}{2}\sqrt{2-\underbrace{\sqrt{2+\sqrt{2+\sqrt{2+\cdots +\sqrt{2}}}}}_{k-1\,\,\text{square}\,\,\text{roots}}}.
\end{equation}
Consequently, from the equations \eqref{eq_6}, \eqref{eq_7} and \eqref{eq_9} we obtain a simple formula for the constant pi
\[
\frac{\pi }{{{2}^{k+1}}}=\arctan \left( \frac{\sqrt{2-\underbrace{\sqrt{2+\sqrt{2+\sqrt{2+\cdots +\sqrt{2}}}}}_{k-1\,\,\text{square}\,\,\text{roots}}}}{\underbrace{\sqrt{2+\sqrt{2+\sqrt{2+\cdots +\sqrt{2}}}}}_{k\,\,\text{square}\,\,\text{roots}}} \right)
\]
or
\begin{equation}\label{eq_10}
\frac{\pi }{{{2}^{k+1}}}=\arctan \left( \frac{\sqrt{2-{{a}_{k-1}}}}{{{a}_{k}}} \right).
\end{equation}
Lastly, combining equations \eqref{eq_5} and \eqref{eq_10} together results in the formula \eqref{eq_2} for pi.

\section{Algorithmic implementation}

\subsection{Methodology description}

As it has been reported previously in the paper \cite{Abrarov2016b}, the decrease of the argument $x$ in the limit \eqref{eq_4a} improves significantly the accuracy in computing pi. Therefore, we may also expect a considerable improvement in accuracy of the arctangent function identity \eqref{eq_5} when its argument $x$ decreases. In fact, the equation \eqref{eq_2} is based on the arctangent function identity \eqref{eq_5} when its argument $x$ is equal to $\sqrt{2-{{a}_{k-1}}}/{{a}_{k}}$. The increment of the integer $k$ by one decreases the argument $x=\sqrt{2-{{a}_{k-1}}}/{{a}_{k}}$ by a factor that tends to two as $k \to \infty$. Therefore, the value of argument $x=\sqrt{2-{{a}_{k-1}}}/{{a}_{k}}$ decreases very rapidly in a geometric progression as the integer $k$ increases. As a consequence, this approach leads to a significant improvement in accuracy of the constant pi.

\subsection{Computational results}

In order to estimate the convergence rate, we performed sample computations of the constant pi by using a rearranged form of the equation \eqref{eq_2} as follows
\footnotesize
\begin{equation}\label{eq_11}
\begin{aligned}
&\hspace{-0.2cm}\pi = \\
&\hspace{+0.1cm}{{2}^{k+1}}i\sum\limits_{m=1}^{m_{max}}{\frac{1}{2m-1}\left( \frac{1}{{{\left( 1+2i\,\,{{a}_{k}}/\sqrt{2-{{a}_{k-1}}} \right)}^{2m-1}}}-\frac{1}{{{\left( 1-2i\,\,{{a}_{k}}/\sqrt{2-{{a}_{k-1}}} \right)}^{2m-1}}} \right)}\\
&+\varepsilon,
\end{aligned}
\end{equation}
\normalsize
where $m_{max} >> 1$ is the truncating integer and $\varepsilon $ is the error term. 

The computational test reveals that even at smallest values of the integer $k$ the equation \eqref{eq_11} can be quite rapid in convergence. In particular, at $m_{max} = 51$ and $k$ equals only to $1$, $2$ and $3$ we can observe a relatively large overlap in digits coinciding with actual value of the constant pi as given by (for $k = 1$ we imply that $a_0 = 0$)
$$
\underbrace{3.141592653\ldots 2795028841}_{38\,\,\text{coinciding}\,\text{digits}}3610032370\ldots
$$
$$
\underbrace{3.141592653\ldots 9230781640}_{73\,\,\text{coinciding}\,\,\text{digits}}5927185386\ldots
$$
and
$$
\underbrace{3.141592653\ldots 1706798214}_{105\,\,\text{coinciding}\,\,\text{digits}}0598570306\ldots \, ,
$$
respectively.

Further, in algorithmic implementation we incremented $k$ by one at each successive step while keeping value of the parameter $m_{max}$ fixed and equal to $169$. Thus, the computational test we performed shows that at $k$ equal to $15$, $16$, $17$, $18$, $19$ and $20$, the error term $\varepsilon $ by absolute value becomes equal to $1.92858\times {{10}^{-1564}}$, $3.44437\times {{10}^{-1666}}$, $6.15140\times {{10}^{-1768}}$, $1.098591\times {{10}^{-1869}}$, $1.96199\times {{10}^{-1971}}$ and $3.50396\times {{10}^{-2073}}$, respectively. As we can see, with only $m_{max}=169$ summation terms each increment of the integer $k$ just by one contributes for more than $100$ additional decimal digits of pi. 

Since the convergence improves while $\sqrt{2-{{a}_{k-1}}}/{{a}_{k}}$ decreases, it is not necessary to increment continuously the integer $k$. For example, when $k$ is fixed and equal to $23$, $50$ and $90$, each increment of the integer $m_{max}$ by one contributes for $14$, $30$ and $54$ decimal digits of pi, respectively.

\section{Theoretical analysis}

\subsubsection*{Proposition}

\noindent According to computational test that has been shown in the previous section, even if the parameter $m_{max}$ is fixed at $169$, the accuracy of pi, nevertheless, improves continuously while the integer $k$ increases. Therefore, relying on these experimental results we assume that the equation \eqref{eq_2} can be modified as
\small
\begin{equation}\label{eq_12}
\begin{aligned}
\pi =& \,\, i\underset{k\to \infty }{\mathop{\lim }}\,{{2}^{k+1}}\sum\limits_{m=1}^{m_{max}}\frac{1}{2m-1}\\
&\times\left( \frac{1}{{{\left( 1+2i\,\,{{a}_{k}}/\sqrt{2-{{a}_{k-1}}} \right)}^{2m-1}}}-\frac{1}{{{\left( 1-2i\,\,{{a}_{K}}/\sqrt{2-{{a}_{k-1}}} \right)}^{2m-1}}} \right).
\end{aligned}
\end{equation}
\normalsize

\subsubsection*{Proof}

The proof is not difficult. Consider the following integral
\footnotesize
\[
\int\limits_{0}^{1}2^{k+1}{\frac{\sqrt{2-{{a}_{k-1}}}/{{a}_{k}}}{1-{{\left( \sqrt{2-{{a}_{k-1}}}/{{a}_{k}} \right)}^{2}}{{t}^{2}}}dt}=2^{k+1}\arctan \left( \frac{\sqrt{2-{{a}_{k-1}}}}{{a}_{k}} \right)=\pi.
\]
\normalsize
The integrand of this integral
\[
g_{k}\left(t\right)=2^{k+1}{\frac{\sqrt{2-{{a}_{k-1}}}/{{a}_{k}}}{1-{{\left( \sqrt{2-{{a}_{k-1}}}/{{a}_{k}} \right)}^{2}}{{t}^{2}}}}
\]
at the limit when $k \to \infty$ becomes
$$
\underset{k\to \infty }{\mathop{\lim }} g_{k}\left(t\right)=g_{\infty}\left(t\right)=\pi.
$$

Since the function $g_{\infty}\left(t\right)$ is just a constant, only its zeroth order of the derivative $g_{\infty}^{\left(0\right)}\left(t\right)$ is not equal to zero. This signifies that if the function $g_{\infty}\left(t\right)$ is substituted into equation \eqref{eq_3b}, then it is no longer necessary to tend the integer $M$ to infinity because for any other than zeroth order of the derivative we always get 
$$
g_{\infty}^{\left(m\right)}\left(t\right) = 0, \qquad m>0.
$$
Consequently, we can infer
\begin{equation}\label{eq_13}
\begin{aligned}
\pi &=\underset{k\to \infty }{\mathop{\lim }}\int\limits_{0}^{1}{g_{k}\left( t \right)dt}\\
&=\underset{M\to \infty }{\mathop{\lim }} \, \underset{k\to \infty }{\mathop{\lim }} \,\sum\limits_{\ell =1}^{L}{\sum\limits_{m=0}^{M}{\frac{{{\left( -1 \right)}^{m}}+1}{{{\left( 2L \right)}^{m+1}}\left( m+1 \right)!}}}{{\left. {{g_{k}}^{\left( m \right)}}\left( t \right) \right|}_{t=\frac{\ell -1/2}{L}}}\\
&=\underset{k\to \infty }{\mathop{\lim }} \,\sum\limits_{\ell =1}^{L}{\sum\limits_{m=0}^{M}{\frac{{{\left( -1 \right)}^{m}}+1}{{{\left( 2L \right)}^{m+1}}\left( m+1 \right)!}}}{{\left. {{g_{k}}^{\left( m \right)}}\left( t \right) \right|}_{t=\frac{\ell -1/2}{L}}}
\end{aligned}
\end{equation}
and from the equations \eqref{eq_13} and \eqref{eq_3b} it follows now that
\footnotesize
\[
\tag{14}\label{eq_14}
\hspace{-8cm}\pi=\, i\underset{k\to \infty }{\mathop{\lim }}2^{k+1}\,\sum\limits_{\ell =1}^{L}\sum\limits_{m=1}^{\left\lfloor \frac{M}{2} \right\rfloor +1}\frac{1}{2m-1}
\]
\[
\hspace{0.7cm}\times\left( \frac{1}{{{\left( \left( 2\ell -1 \right)+2iL\,\,{{a}_{k}}/\sqrt{2-{{a}_{k-1}}} \right)}^{2m-1}}}-\frac{1}{{{\left( \left( 2\ell -1 \right)-2iL\,\,{{a}_{k}}/\sqrt{2-{{a}_{k-1}}} \right)}^{2m-1}}} \right).
\]
\normalsize
\\
\noindent 
Implying now that $m_{max}=\left\lfloor M/2 \right\rfloor +1$, at $L=1$ this limit is reduced to equation \eqref{eq_12}. This completes the proof.

\subsubsection*{Corollary}

At $L=1$ and $M=1$ the limit \eqref{eq_14} is simplified as
\setcounter{equation}{14}
\small
\begin{equation}\label{eq_15}
\begin{aligned}
\pi = \,\, i\underset{k\to \infty }{\mathop{\lim }}\,2^{k+1}\left( \frac{1}{{{ 1+2i\,\,{{a}_{k}}/\sqrt{2-{{a}_{k-1}}}}}}-\frac{1}{{{1-2i\,\,{{a}_{k}}/\sqrt{2-{{a}_{k-1}}}}}} \right)
\end{aligned}
\end{equation}
\normalsize
and since (see equation \eqref{eq_6})
\[
\underset{k\to \infty }{\mathop{\lim }} a_k = 2 \underset{k\to \infty }{\mathop{\lim }} \cos \left( \frac{\pi}{2^{k+1}} \right) = 2
\]
the limit \eqref{eq_15} provides
\small
\begin{equation}\label{eq_16}
\begin{aligned}
\pi =& \,\, i\underset{k\to \infty }{\mathop{\lim }}\,2^{k+1}\left( \frac{1}{{{ 1+4i\,\,/\sqrt{2-{{a}_{k-1}}}}}}-\frac{1}{{{1-4i\,\,/\sqrt{2-{{a}_{k-1}}}}}} \right)\\
=& \underset{k\to \infty }{\mathop{\lim }}\,2^{k+1} \frac{8/\sqrt{2-{{a}_{k-1}}}}{1+16 /\left({2-{{a}_{k-1}}}\right)}\\
=& \underset{k\to \infty }{\mathop{\lim }}\,2^{k+1} \frac{8}{\sqrt{2-{{a}_{k-1}}}+16 /\sqrt{2-{{a}_{k-1}}}}.
\end{aligned}
\end{equation}
\normalsize
From the equation \eqref{eq_8} it immediately follows that
\[
\underset{k\to \infty }{\mathop{\lim }} \sqrt{2-a_{k-1}} =2 \underset{k\to \infty }{\mathop{\lim }} \sin \left( \frac{\pi}{2^{k+1}} \right) = 0.
\]
Consequently, the limit \eqref{eq_16} can be simplified further as
\[
\pi = \underset{k\to \infty }{\mathop{\lim }}\,2^{k+1} \frac{8}{16 /\sqrt{2-{{a}_{k-1}}}}=\underset{k\to \infty }{\mathop{\lim }}\,2^{k} \sqrt{2-a_{k-1}}
\]
or
\[
\pi =\underset{k\to \infty }{\mathop{\lim }}\,{{2}^{k}}\sqrt{2-\underbrace{\sqrt{2+\sqrt{2+\sqrt{2+\cdots +\sqrt{2}}}}}_{k-1\,\,\text{square}\,\,\text{roots}}}.
\]
This is a well-known formula for pi \cite{Servi2003}.

\section{Conclusion}

A new formula \eqref{eq_2} for pi based on the arctangent function identity \eqref{eq_5} with argument $x=\sqrt{2-{{a}_{k-1}}}/{{a}_{k}}$, where
\[
{{a}_{k}}=\underbrace{\sqrt{2+\sqrt{2+\sqrt{2+\cdots +\sqrt{2}}}}}_{k\,\,\text{square}\,\,\text{roots}},
\]
is presented. This approach demonstrates high efficiency in computation due to rapid convergence. Specifically, the computational test reveals that with only $169$ summation terms the increment of integer $k$ just by one provides more than $100$ additional decimal digits of the constant pi.

\section*{Acknowledgments}

This work is supported by National Research Council Canada, Thoth Technology Inc. and York University. The authors thank the reviewers for constructive comments and recommendations.

\bigskip



\begin{thebibliography}{9}

\bibitem{Servi2003}
L.D. Servi, Nested square roots of 2, Amer. Math. Monthly, 110 (4) (2003) 326-330. \\
\url{http://dx.doi.org/10.2307/3647881}

\bibitem{Levin2005}
A. Levin, A new class of infinite products generalizing Vi\'{e}te\text{'}s product formula for $\pi$, Ramanujan J. 10 (3) (2005) 305-324. \\
\url{http://dx.doi.org/10.1007/s11139-005-4852-z}

\bibitem{Kreminski2008}
R. Kreminski, $\pi$ to thousands of digits from Vieta\text{'}s formula, Math. Magazine, 81 (3) (2008) 201-207. \\
\url{http://www.jstor.org/stable/27643107}

\bibitem{Abrarov2016a}
S.M. Abrarov and B.M. Quine, Identities for the arctangent function by enhanced midpoint integration and the high-accuracy computation of pi, \href{http://arxiv.org/abs/1604.03752}{arXiv:1604.03752}, 2016.

\bibitem{Abrarov2016b}
S.M. Abrarov and B.M. Quine, A simple identity for derivatives of the arctangent function, \href{http://arxiv.org/abs/1605.02843}{arXiv:1605.02843}, 2016.

\end{thebibliography}
\end{document}